\documentclass[11pt]{article}

\input{amssym.def}
\input{amssym.tex}

\title{Simon's conjecture for 2-bridge knots}
\author{Michel Boileau\thanks{supported in part by the ANR},~~Steve Boyer\thanks{supported in part by NSERC},~~Alan W. Reid\thanks{supported in part
by the NSF},~~Shicheng Wang\thanks{supported in part by grant no.
10631060 of NSF of China.}}

\def\ker{\rm{ker}}

\def\dim{\rm{dim}}

\def\PSL{\rm{PSL}}
\def\SU{\rm{SU}}
\def\SL{\rm{SL}}

\def\qed{ $\sqcup\!\!\!\!\sqcap$}

\def\tr{\mbox{\rm{tr}}}

\newtheorem{theorem}{Theorem}[section]
\newtheorem{lemma}[theorem]{Lemma}
\newtheorem{corollary}[theorem]{Corollary}
\newtheorem{proposition}[theorem]{Proposition}
\newtheorem{conjecture}[theorem]{Conjecture}
\newtheorem{question}[theorem]{Question}
\newtheorem{definition}[theorem]{Definition}

\begin{document}

\maketitle \begin{abstract} It is conjectured that for each knot $K$
in $S^3$, the fundamental group of its complement
surjects onto only finitely many distinct knot groups.
Applying character variety theory we obtain
an affirmative solution of
the conjecture for a class of small knots that includes
2-bridge knots.\end{abstract}

%
%
%
%
\section{Introduction}

In this paper all knots and links are in $S^3$. For a knot (link)
$K$, we often simply call the fundamental group of $S^3\setminus K$,
{\em the group of $K$} or {\em the knot (link) group.}
Let $K$ be a non-trivial knot. Simon's Conjecture (see \cite{Ki}
Problem 1.12(D)) asserts the following:

\begin{conjecture}
\label{simons}
$\pi_1(S^3\setminus K)$ surjects onto only finitely many distinct knot groups.\end{conjecture}

Although this conjecture dates back to the $1970$'s, and has received
considerable attention recently (see \cite{BRW}, \cite{BB},
\cite{GAR3}, \cite{ORS}, \cite{RW}, \cite{RWZ}, \cite{SW} and
\cite{SW1} to name a few), little by way of general results appears to
be known. Conjecture \ref{simons} is easily seen to hold for torus
knots (we give the proof in \S 3.1). In \cite{BRW}, the conjecture is
established under the assumption that the epimorphisms are
non-degenerate in the sense that the longitude of $K$ is sent to a
non-trivial peripheral element under the epimorphism. In particular
this holds in the case when the homomorphism is induced by mapping of
non-zero degree.

Since any knot group is the homomorphic image of
the group of a hyperbolic knot (see for example \cite{Ka}), it is
sufficient to prove the conjecture for hyperbolic knots.
The main result of this paper is the following, and is the first
general result for a large class of hyperbolic knots.

\begin{theorem}
\label{main}
Conjecture \ref{simons} holds for all 2-bridge knots.\end{theorem}

Indeed, using \cite{BRW},  \cite{BB}, \cite{GAR1} and
\cite{GAR3},
one can say more about the nature of the homomorphisms in Theorem \ref{main}.
In particular, we establish:

\begin{corollary}
\label{nonzerodegree} Let $K$ be a 2-bridge hyperbolic knot, and $K'$
a non-trivial knot. If there is an epimorphism
$\varphi:\pi_1(S^3\setminus K)\rightarrow \pi_1(S^3\setminus K')$,
then $\varphi$ is induced by a map  $f: S^3\setminus K \rightarrow
S^3\setminus K'$ of non-zero degree. Furthermore, $K'$ is
necessarily a 2-bridge knot.
\end{corollary}

As we discuss in \S 4.2, a strengthening of Corollary \ref{nonzerodegree}
holds, where, epimorphism is replaced by {\em a virtual epimorphism} (see
\S 4.2 for more details).

We will also prove some general results towards Conjecture \ref{simons} for a
larger class of knots satisfying certain conditions. These results will
be used in proving Theorem \ref{main} and Corollary
\ref{nonzerodegree}. Recall that a compact oriented 3-manifold $N$
is called {\it small} if $N$ contains no closed embedded essential
surface; and a knot (link) $K\subset S^3$ is called {\em small} if
the exterior $E(K)$ is small. The fact that 2-bridge knots are small is proved in
\cite{HT}.

\begin{theorem}
\label{smallsimonintro} Let $L$ be a small hyperbolic link of n
components. Then $\pi_1(S^3\setminus L)$ surjects onto only finitely
many groups of hyperbolic links of $n$ components.\end{theorem}

Note that a knot is either a torus knot, or a hyperbolic knot, or a
satellite knot.
In particular, it follows from Theorem \ref{smallsimonintro}, that
if $K$ is a small hyperbolic knot, then
$\pi_1(S^3\setminus K)$ surjects onto only finitely
many hyperbolic knot groups. As we discuss below, it is easy to
establish using the Alexander polynomial, that any knot group
surjects onto only finitely many distinct torus knot groups.

It is perhaps tempting at this point to think that there cannot be an
epimorphism from a small knot group to the group of a satellite
knot, therefore Conjecture \ref{simons} holds for small knots.  This does
not seem so easy to exclude and motivates the following more general
question:

\begin{question}
\label{small_surject_large}
Does there exist a small knot  $K\subset S^3$ such that $\pi_1(S^3\setminus K)$
surjects onto $\pi_1(S^3\setminus K')$, and $S^3\setminus K'$ contains a closed embedded essential
surface?\end{question}

When the target is the fundamental group of a satellite knot, we
will prove that if such a homomorphism exists, then the longitude
$\lambda$ of $K$ must be in the kernel. More precisely,

\begin{proposition}
\label{longnotdead} Let $K$ be a small hyperbolic knot and let $K'$
be a satellite knot. Assume that $\varphi:\pi_1(S^3\setminus
K)\rightarrow \pi_1(S^3\setminus K')$ is an epimorphism. Then
$\varphi(\lambda)=1$.\end{proposition}

In order to control the image of the longitude, we introduce the
following definition which can be thought as a kind of smallness for
the knot in terms of the character variety of longitudinal surgery
on $K$. For a knot $K\subset S^3$
we denote by $K(0)$ the manifold obtained from $S^3$ by a longitudinal
surgery on $K$.

\begin{definition}
\label{propertyL}
Let $K$ be a knot. We will say $K$ has Property L if the
$\SL(2,{\bf C})$-character variety of the manifold $K(0)$ contains
only finitely many characters of irreducible representations.
\end{definition}

The motivation for this definition is the following result.

\begin{proposition}
\label{nontrivlong}
If a hyperbolic knot has Property L, then for any non-trivial knot $K'$ and epimomorphism
$\varphi:\pi_1(S^3\setminus K)\rightarrow \pi_1(S^3\setminus K')$, $\ker~\varphi$ does not contain the longitude of $K$.
\end{proposition}

The following result is thus a
consequence of Theorem \ref{smallsimonintro} and Propositions
\ref{longnotdead} and \ref{nontrivlong}:

\begin{theorem}
\label{mainsmall}
Let $K$ be a small knot and assume that $K$ has Property L.
Then Conjecture \ref{simons} holds
for $K$.\end{theorem}

Our main result (Theorem \ref{main}) now follows immediately
from Theorem \ref{mainsmall} and the next proposition.

\begin{proposition}
\label{2bridgehasL} Let $K$ be a hyperbolic 2-bridge knot. Then $K$
has Property L.\end{proposition}

Although Property L is framed in terms of the character variety, which
can be difficult to understand, there are useful criteria which are
sufficient for a small knot to have Property L. The first one will be
used to show that Property L holds for 2-bridge knots. For the
definition of a parabolic representation or of a strict boundary slope
see \S 2 .

\begin{proposition}
\label{criterpropertyL}
Let $K$ be a small hyperbolic knot. \begin{itemize}
\item If no parabolic representation $\rho: \pi_1(S^3\setminus K) \to
\SL(2,{\bf C})$ kills the longitude of $K$, then property $L$ holds for $K$.
\item If the longitude is not a strict boundary slope, then
Property L  holds for $K$.\end{itemize}
\end{proposition}

\noindent{\bf Remarks on Propery L:} (1) Using Proposition
\ref{nontrivlong}, it is easy to construct knots
which do not have Property L.

For example, using the construction of
\cite{Joh} on a normal generator for a knot group that is not a
meridian (which exist in some abundance \cite{BRW}, see \cite{Ca} for
explicit examples), one can construct a hyperbolic knot whose group
surjects onto another hyperbolic knot group sending the longitude
trivially.  In \cite{GAR3}, examples are given where the domain knot
is small; for example there is an epimorphism of the group of the knot
$8_{20}$ onto the group of the trefoil-knot for which the longitude of
$8_{20}$ is mapped trivially.\\[\baselineskip]
(2) Control of the image of the longitude has featured in other work
related to epimorphisms between knot groups;  for example Property
$Q^*$ of Simon (see \cite{Si} and also \cite{GAR1}, \cite{GAR3}).
Indeed, from \cite{GAR3}, the property given by Proposition \ref{nontrivlong} can be viewed as an extension of Property $Q^*$ of Simon.\\[\baselineskip]
(3) Note that if $K$ and $K'$ are knots with Alexander polynomials
$\Delta_K(t)$ and $\Delta_{K'}(t)$ respectively, and
$\varphi:\pi_1(S^3\setminus K)\rightarrow \pi_1(S^3\setminus K')$ an
epimorphism, then it is well-known that
$\Delta_{K'}(t)|\Delta_K(t)$. Thus, simple Alexander polynomial
considerations shows that any knot group surjects onto only finitely
many distinct torus knot groups, and so it is only when the target
is hyperbolic or satellite that the assumption of Property L is
interesting.\\[\baselineskip]
The character variety (as in \cite{BB} and \cite{RW}) is the main
algebraic tool that organizes the proofs of the results in this paper.
In particular we make use of the result of Kronheimer and Mrowka
\cite{KM} which ensures that the $\SL(2,{\bf C})$-character variety
(and hence the $\PSL(2,{\bf C})$-character variety)
of any non-trivial
knot contains a curve of characters of irreducible
representations.\\[\baselineskip]
\noindent{\bf A comment on application of the character variety to Simon's
conjecture:} As we can and will see  from \cite{BB}, \cite{RW} and
the present paper, the theory of character varieties is particularly
useful in the study of epimorphisms between 3-manifolds groups when the
domain manifolds are small.

However, comparison of the two results below suggests possible limitation of
applying character variety methods to Simon's conjecture as well as the
truth of Simon's conjecture itself: On the one hand the group of each
small hyperbolic link of n components surjects onto only finitely
many groups of $n$ component hyperbolic links  (Corollary
\ref{smallsimon}); while on the other hand, there exist hyperbolic
links of two components whose groups surject onto
the group of every two bridge link
(see the discussion related to Conjecture
\ref{nSimon}).\\

\noindent{\bf Organization of the paper:} The facts about the
character variety that will be used later are presented in Section 2. Results
stated for small knots, such as Theorem \ref{smallsimonintro},
Proposition \ref{longnotdead}, Proposition \ref{nontrivlong},
Theorem \ref{mainsmall} and Corollary
\ref{notstrictboundaryslope}, will be proved in Section 3. Results
stated for 2-bridge knots, such as Theorem \ref{main}, Proposition
\ref{nonzerodegree}, Proposition \ref{2bridgehasL},  will be proved
in Section 4. Section 5 records more questions, consequences and
facts for the character variety and Simon's conjecture that have
arisen out
of our work.\\[\baselineskip]
\noindent{\bf Acknowledgements:}~The second and third authors wish
to thank the Department of Mathematics at  Universit\'e Paul
Sabatier for their hospitality during this work. The third author
also wishes to thank the Department of Mathematics at Peking
University where this work started, and The Institute for Advanced Study
where this work continued.

\section{Preliminaries}

\subsection{Some notation}

\noindent Throughout, if $L\subset S^3$ is a link we shall let
$E(L)$ denote the exterior of $L$; that is the closure of the
complement of a small open tubular neighbourhood of $L$. If
$K\subset S^3$ is a knot and $r\in{\bf Q}\cup \infty$ a slope, then
$K(r)$ will denote the manifold obtained by $r$-Dehn surgery on $K$
(or equivalently, $r$-Dehn filling on $E(K)$). Our convention is
always that a meridian of $K$ has slope $1/0$ and a longitude $0/1$.
A slope $r$ is called a {\em boundary slope}, if $E(K)$ contains an
embedded essential surface whose boundary consists of a non-empty
collection of parallel copies of simple closed curves on $\partial
E(K)$ of slope $r$.  The longitude of a knot $K$ always bounds
Seifert surface of $K$, and so is a boundary slope. It is called a
{\em strict boundary slope} if it is the boundary
slope of a surface that is not a fiber in a fibration over the circle.

\subsection{Standard facts about the character variety}

 Let $G$ be a finitely generated group. We denote by $X(G)$ (resp. $Y(G)$) the
$\SL(2,{\bf C})$-character variety (resp. $\PSL(2,{\bf
C})$-character variety) of $G$ (see \cite{CS} and \cite{BoZ1} for
details). If $V$ is an algebraic set, we define the {\em dimension}
of $V$ to be the maximal dimension of an irreducible component of
$V$. We will denote this by $\dim (V)$.

 Suppose that $G$ and $H$ are finitely generated groups and
 $\varphi:G\rightarrow H$ is an epimorphism.  Then $\varphi$ defines a map at the level of character
 varieties $\varphi^*:X(H)\rightarrow X(G)$ by $\varphi^*(\chi_\rho) = \chi_{\rho\circ \varphi}$. This
 map is algebraic, and furthermore is a closed map in the Zariski topology (see \cite{BB} Lemma 2.1).
 In future we will abbreviate composition of homomorphisms $\varphi\circ\psi$ by $\varphi\psi$.

We make repeated use of the following standard fact.

  \begin{lemma}
 \label{one-to-one}
 Let $G$ and $H$ be as above,  then $\varphi^*$ injects $X(H)\hookrightarrow X(G)$.\end{lemma}

 \noindent{\bf Proof:}~Suppose $\chi_\rho,\chi_{\rho'}\in X(H)$ with  $\varphi^*(\chi_\rho)=\varphi^*(\chi_{\rho'})$.
 Thus, $\chi_{\rho\varphi}(g) = \chi_{\rho'\varphi}(g)$ for all $g\in G$, and since $\varphi$ is onto, we deduce
 that $\chi_\rho(h) = \chi_{\rho'}(h)$ for all $h\in H$. Hence $\chi_\rho = \chi_{\rho'}$.\qed\\[\baselineskip]
 We now assume that $D\subset X(G)$ is a component containing the
character $\chi_\rho$ of an irreducible representation and
$D=\varphi^*(C)$ (as noted $\varphi^*$ is a closed map) for some
component $C\subset X(H)$. Then, $\chi_\rho =
\varphi^*(\chi_{\rho'})$ for some irreducible representation
$\rho':H\rightarrow \SL(2,{\bf C})$. By definition,
$\varphi^*(\chi_{\rho'})=\chi_{\rho'\varphi}$, and so since the
representations $\rho$ and $\rho'\varphi$ are irreducible, we deduce
that  the groups $\rho(G)$ and $\rho'\varphi(G) = \rho'(H)$ are
conjugate in $\SL(2,{\bf C})$. In particular, after conjugating if
necessary, the homomorphisms $\rho' \varphi$ and $\rho$ have the
same image.

\subsection{Existence of irreducible representations of knot groups}

When $G=\pi_1(M)$, and $M$ is a compact 3-manifold we denote $X(G)$
(resp. $Y(G)$) by $X(M)$ (resp. $Y(M)$). When $M$ is a knot exterior
in $S^3$ we write $X(M)=X(K)$ (resp. $Y(M)=Y(K)$).

Now $X(K)$ (resp. $Y(K)$) always contains a curve of characters
corresponding to abelian representations.  When $K$ is a hyperbolic
knot (i.e. $S^3\setminus K$ admits a complete hyperbolic structure of
finite volume), it is a well-known consequence of Thurston's Dehn
surgery theorem (see \cite[Proposition 1.1.1]{CGLS}) that there is a
so-called {\em canonical component} in $X(K)$ (resp. $Y(K)$) which
contains the character of a faithful discrete representation of
$\pi_1(S^3\setminus K)$.  More recently, the work of Kronheimer and
Mrowka \cite{KM} establishes the following general result (we include
a proof of the mild extension of their work that is needed for us).

 \begin{theorem}
 \label{KM}
 Let $K$ be a non-trivial knot. Then $X(K)$ (resp. $Y(K)$) contains a
 curve for which all but finitely many of its elements are characters
 of irreducible representations.
 \end{theorem}

\noindent{\bf Proof:}~It suffices to prove Theorem \ref{KM} for $X(K)$. As
the set of reducible characters is Zariski closed in $X(K)$
(\cite[proof of Corollary 1.4.5]{CS}), by a result of Thurston (see
\cite[Proposition 3.2.1]{CS}) to find a curve in the conclusion of
Theorem \ref{KM}, it is enough to find an irreducible representation
$\rho: \pi_1(E(K)) \to \SL(2;{\bf C})$ such that
$\rho(\pi_1(\partial E(K))) \not \subset \{\pm I\}$. Note that the
latter condition holds for any irreducible representation of
$\pi_1(E(K))$.

To find an irreducible $\rho$, note that by \cite{KM}, for any
$r\in{\bf Q}$ with $|r| \leq 2$, $\pi_1(K(r))$, admits a non-cyclic
$\SU(2)$-representation. Take $r=1$ and suppose that the
representation guaranteed by \cite{KM} is reducible as a
representation into $\SL(2,{\bf C})$. Since $\pi_1(K(1))$ is
perfect, it coincides with its commutator subgroup and therefore the
trace of any element of the image of $\rho$ is $2$. As $I$ is the
only element of $\SU(2)$ with this trace, the image of $\rho$ is
$\{I\}$, a contradiction. \qed\\


\subsection{$X(K)$ for small hyperbolic knots and p-rep. characters}

We now prove some results about the character variety of a small
hyperbolic knot.  It will be convenient to recall some terminology
from \cite{CGLS}.

Let $K\subset S^3$ be a knot and $\alpha\in \pi_1(\partial E(K))$.  If
$X\subset X(K)$ is a component, define the polynomial function:

$$f_\alpha:X\rightarrow {\bf C}~~\hbox{by}~~f_\alpha(\chi_\rho) = \tr^2(\rho(\alpha))-4.$$

We first record the following well-known result.

\begin{theorem}
\label{move_meridian} (1) Let $N$ be hyperbolic 3-manifold with
$\partial N$ a union of $n$ tori. If $X$ is an irreducible component
of $X(N)$ that contains the character of an irreducible
representation, then $\dim (X)$ is at least $n$; moreover $\dim
(X)=n$ when $N$ is small.

(2) Let $K$ be a small hyperbolic knot and $\mu$ be a meridian of
$K$. If $x$ is an ideal point of $X$, then $f_\mu$ has a pole at $x$. In
particular, $f_\mu$ is non-constant.\end{theorem}

\noindent{\bf Proof:}~ (1) The dimension of $X$ is at least $n$ by
\cite[Proposition 3.2.1]{CS} and at most $n$ when $N$ is small by
\cite[Theorem 4.1]{CL}.

(2) Let $x$ be an ideal point of $X$ and consider
$$I_\mu : X \rightarrow {\bf C}, \;\;\; I_\mu(\chi_\rho)=\tr(\rho(\mu)).$$
Clearly $f_\mu = I_\mu^2 - 4$, so to prove the lemma it suffices to
show that $I_\mu$ has a pole at $x$. Now  \cite[Proposition
1.3.9]{CGLS} implies that either $I_\mu(x) = \infty$, or $\mu$ is a
boundary slope, or $I_\alpha(x) \in {\bf C}$ for all $\alpha \in
\pi_1(\partial E(K))$. The second possibility is ruled out by
\cite[Theorem 2.0.3]{CGLS}, while the third is ruled out by the fact
that it implies $E(K)$ contains a closed essential surface (cf. the
second paragraph of \cite[\S 1.6.2]{CGLS}),
which contradicts that $E(K)$ is small.\qed\\

\noindent Note that zeroes of $f_\alpha$ correspond to
representations $\rho$ for which $\alpha$ either maps trivially (in
$\PSL(2,{\bf C})$) or to a parabolic element.  In this latter case,
it is easy to see that $f_\beta(\chi_\rho)=0$ for all
$\beta\in\pi_1(\partial E(K))$. Following Riley \cite{Ri1}, we call
such a representation a  {\em parabolic representation} or {\em
p-rep}.  We define a character $\chi_\rho$ to be a {\em p-rep
character} if $\rho$ is an irreducible representation for which at
least one peripheral element is mapped to a parabolic element.\\

The following proposition will be useful.

\begin{proposition}
\label{prepchars}
Let $K$ be a small hyperbolic knot and $X\subset X(K)$ an irreducible
component that contains the character of an irreducible
representation.  Then $X$ contains a p-rep character. Indeed, the set
of $p$-rep characters on $X$ is the zero set of $f_\mu$ on $X$.\end{proposition}

\noindent{\bf Proof:}~By Theorem \ref{move_meridian} (1), $X$ is a
curve. Let $\tilde{X}$ be its smooth projective model. Then
$\tilde{X} = X^\nu \cup \mathcal{I}$ where $\nu: X^\nu \to X$ is an
affine desingularisation and $\mathcal{I}$ is the finite set of
ideal points of $X$. The function $f_\mu$ corresponds to a
holomorphic map $\tilde{f_\mu}:\tilde{X}\rightarrow {\bf CP}^1$ (see
\cite{CS}) where $\tilde{f_\mu}|X^\nu = f_\mu \circ \nu$. Thus
Theorem  \ref{move_meridian} implies that $\tilde{f_\mu}$ is
non-constant, so it has at least one zero $x_0$, and also that $x_0
\in X^\nu$. Set $\nu(x_0) = \chi_\rho$. Since $X$ contains an
irreducible character, \cite[Proposition 1.5.5]{CGLS} implies that
we can suppose the image of $\rho$ is non-cyclic. Hence $\rho(\mu)
\ne \pm I$ and therefore $\rho(\mu)$ is parabolic. It follows that
if $\alpha \in \pi_1(\partial E(K))$, then either $\rho(\alpha)$ is parabolic,
or $\rho(\alpha)$ is $\pm I$. Thus, the proof of the proposition will be
complete once we establish that $\rho$ is irreducible.

Suppose this were not the case and let $R$ be the $4$-dimensional
component of the representation variety $R(K) =
\hbox{Hom}(\pi_1(E(K)), \SL(2; {\bf C}))$ whose image in $X(K)$ equals
$X$ (cf. \cite[Lemma 4.1]{BoZ1} and \cite[Proposition 1.5.3]{CS}).  By
\cite[Proposition 1.5.6]{CGLS} we can suppose $\rho \in R$.  Since $R$
is $\SL(2;{\bf C})$-invariant (\cite[Proposition 1.1.1]{CS}), we can
suppose that the image of $\rho$ consists of upper-triangular
matrices.  Hence consideration of the sequence $\rho_n =
\left( \begin{array}{cc} \frac{1}{n} & 0 \\ 0 & n \end{array} \right)
\rho \left( \begin{array}{cc} \frac{1}{n} & 0 \\ 0 & n \end{array}
\right)^{-1}$ shows that $R$ contains a representation $\rho_0$ whose
image is diagonal and which sends $\mu$ to $\pm I$. Thus
$\rho_0(\gamma) = \pm I$ for all $\gamma \in \pi_1(E(K))$.  The
Zariski tangent space of $R$ at $\rho_0$ is naturally a subspace of
the vector space of $1$-cocycles $Z^1(\pi_1(E(K)); \hbox{sl}(2;{\bf
  C})_{Ad \circ \rho_0})$ (see \cite{We}).  Since the image of
$\rho_0$ is central in $\SL(2; {\bf C})$, $\hbox{sl}(2;{\bf C})_{Ad
  \circ \rho_0}$ is a trivial $\pi_1(E(K))$-module. It follows that
$Z^1(\pi_1(E(K)); \hbox{sl}(2;{\bf C})_{Ad \circ \rho_0}) \cong
H^1(\pi_1(E(K)); \hbox{sl}(2;{\bf C})_{Ad \circ \rho_0}) \cong
H^1(\pi_1(E(K)); {\bf C}^3) \cong {\bf C}^3$.  Hence the dimension of
the Zariski tangent space of $R$ at $\rho_0$ is at most $3$. But this
contradicts the fact that $R$ is $4$-dimensional. Thus $\rho$ must be
irreducible.

To complete the proof, simply note that we have shown that each zero
of $f_\mu$ on a curve component of $X(K)$ containing the character of
an irreducible representation is the character of a p-rep. The
converse is obvious.  \qed

\section{Results for small knots}

\subsection{Simon's Conjecture for torus knots}

In this section we give a quick sketch of the proof that torus knots
satisfy Conjecture \ref{simons} (see also \S 2 of \cite{SW1}). Here
Property L is not needed.

Thus suppose that $K$ is a torus knot, and assume that there exist
infinitely many distinct knots $K_i$ and epimorphisms
$$\varphi_i: \pi_1(S^3\setminus K) \rightarrow  \pi_1(S^3\setminus K_i).$$
Note that if $z$ generates the center of $\pi_1(S^3\setminus K)$, then
$\varphi_i(z)\neq 1$; otherwise, $\varphi_i$ factorizes through a
homomorphism of the base orbifold group $C_{r,s}$ which is the free
product of two cyclic groups of orders $r$ and $s$ for some co-prime
integers $r$ and $s$.  This is impossible, since $\pi_1(S^3\setminus
K_i)$ is torsion-free.  Thus $\pi_1(S^3\setminus K_i)$ has non-trivial
center, and so is a torus knot group by Burde-Zieschang's
characterization of torus knots \cite{BuZ}.

However, as mentioned in \S 1,  $\Delta_{K_i}(t)$ will be a factor of  $\Delta_K(t)$, and so it easily
follows that only finitely many of these $K_i$ can be distinct torus knots.  This completes the proof.\qed

Using this result for torus knots, to prove Conjecture \ref{simons} for
small knots, it therefore suffices to deal with the cases where the domain is a
hyperbolic knot. That is the case we will consider in the remainder
of this section.

\subsection{Proof of Theorem \ref{smallsimonintro}}

In this section we will first prove Theorem \ref{smallsimonintro}. As
remarked upon in \S 1, the finiteness of torus knot groups follows
from Alexander polynomial considerations. The finiteness of
hyperbolic knot group targets follows easily from our next result.
Recall that an
elementary fact in algebraic geometry is that the number of irreducible
components of an algebraic set $V$ is finite, and
hence there are only finitely many of any given dimension $n$.

\begin{theorem}
\label{numberofcptsbdd} Let $G$ be a finitely generated group.
Assume that $\dim (X(G)) = n$ and let $m$ denote the number of
irreducible components of $X(G)$ of dimension $n$.  Suppose that for
$i=1,\ldots ,k$, $N_i$ is a hyperbolic 3-manifold with
incompressible boundary consisting of precisely $n$ torus boundary
components, and that $G$ surjects onto $\pi_1(N_i)$.  We assume that
the $N_i$'s are all non-homeomorphic. Then $k\leq m$.\end{theorem}

\noindent{\bf Proof:}~ Let $\varphi_i:G \rightarrow \pi_1(N_i)$ be the
surjections for
$i=1, \ldots, k$. As discussed in \S 2.1, this induces a closed
algebraic map $\varphi_i^*:X(N_i)\hookrightarrow X(G)$ that is
injective by Lemma \ref{one-to-one}. Furthermore, if $X_i$ denotes
the canonical component of $X(N_i)$, then $\dim (X_i) = n$ by Thurston's
Dehn Surgery Theorem.

Suppose that $k>m$. Then there exists $i,j\in\{1,\ldots, k\}$, and
an irreducible component $X' \subset X(G)$ of dimension at least $n$
such that:
$$\varphi_i^*(X_i), \varphi_j^*(X_j) \subset X'.$$
By the injectivity of $ \varphi_i^*$ and the assumption that $\dim
(X(G)) = n$ it follows that $\varphi_i^*(X_i)$, $\varphi_j^*(X_j)$ and
$X'$ all have dimension $n$ and so
$$\varphi_i^*(X_i)=X'=\varphi_j^*(X_j).$$
Relabelling for convenience, we set $i,j=1,2$.  The equality of
these varieties implies that for each
$$\chi_{\rho_1} \in X_1,~~\hbox{ there exists}~~\chi_{\rho_2'}\in X_2~~\hbox{with}~~
                \varphi_1^*(\chi_{\rho_1}) = \varphi_2^*(\chi_{\rho_2'}),$$
and for each
$$\chi_{\rho_2} \in X_2,~~\hbox{ there exists}~~\chi_{\rho_1'}\in X_1~~\hbox{with}~~
                \varphi_2^*(\chi_{\rho_2}) = \varphi_1^*(\chi_{\rho_1'}).$$
In particular,  we can take $\rho_1$ to be the faithful discrete
representation of $\pi_1(N_1)$, and $\rho_2$ to be the faithful
discrete representation of $\pi_1(N_2)$. Since both $\rho_1$ and
$\rho_2$ are faithful, we have

$$\rho_1(\pi_1(N_1)) \cong \pi_1(N_1)~~\hbox{and}~~\rho_2(\pi_1(N_2)) \cong \pi_1(N_2)).$$
Hence from above, this yields representations
$\rho_1':\pi_1(N_1)\rightarrow \SL(2,{\bf C})$ and
$\rho_2':\pi_1(N_2)\rightarrow \SL(2,{\bf C})$ which satisfy

$$\rho'_2(\pi_1(N_2)) \cong \pi_1(N_1)~~\hbox{and}~~\rho'_1(\pi_1(N_1)) \cong \pi_1(N_2)).$$
Hence, we get {\em epimorphisms}:

$$\rho_2'\rho_1':\pi_1(N_1)\rightarrow \pi_1(N_1)~~\hbox{and}~~\rho_1'\rho_2':\pi_1(N_2)\rightarrow \pi_1(N_2).$$
It is well-known that the fundamental groups of compact hyperbolic
3-manifolds are Hopfian, and so $\rho_2'\rho_1'$ and
$\rho_1'\rho_2'$ are isomorphisms. It now follows that $\rho'_1$
must be also an injection, hence $\pi_1(N_1)\cong \pi_1(N_2)$. Since
both $N_1$ and $N_2$ are complete hyperbolic 3-manifolds with finite
volume,  $N_1$ and $N_2$ are homeomorphic by Mostow Rigidity
Theorem, which contradicts the assumption that they are
non-homeomorphic.\qed\\[\baselineskip]
The most interesting and  immediate application of Theorem
\ref{numberofcptsbdd} is the following:

\begin{corollary}
\label{smallsimon} Let $L$ be a small hyperbolic  link of n
components. Then $\pi_1(S^3\setminus L)$ surjects onto only finitely
many groups of hyperbolic links of $n$ components.\end{corollary}

\noindent{\bf Proof:}~The exterior of each link of
$n$-components has an union of $n$ tori as boundary, and  for a small hyperbolic
link $L$ of $n$ components $\dim (X(L)) = n$ by Theorem
\ref{move_meridian} (1).
Then the proof follows readily from Theorem  \ref{numberofcptsbdd}.\qed\\

\noindent Now Theorem \ref{smallsimonintro} follows from Corollary
\ref{smallsimon} and the discussion about torus knots in
\S 1.\qed\\[\baselineskip]
Theorem \ref{numberofcptsbdd} also provides information about the
nature of $X(K)$ for possible counterexamples to Conjecture
\ref{simons}.

\begin{corollary}
\label{highdimcpt} Suppose $K\subset S^3$ is a hyperbolic knot and
assume that $K_i\subset S^3$ is an infinite family of distinct
hyperbolic knots for which there are epimorphisms
$\varphi_i:\pi_1(S^3\setminus K)\rightarrow \pi_1(S^3\setminus
K_i)$. Then $X(K)$ contains an irreducible component of dimension at
least $2$.\end{corollary}

 \noindent{\bf Proof:}~If  all components
have dimension 1, then Theorem \ref{numberofcptsbdd} bounds the
number of knots $K_i$.\qed\\[\baselineskip]
\noindent{\bf Remark:}~Theorem \ref{numberofcptsbdd} can also be
formulated for the $\PSL(2,{\bf C})$-character variety.

\subsection{Satellite targets}

In this section we prove Proposition \ref{longnotdead}. Before
giving the proof we fix some notation that will be employed in \S 3.3 and
\S 3.4.\\[\baselineskip]
\noindent {\bf Notation:}~Let $K$ be a knot,
$\lambda$ be the longitude for $K$, $\mu$ a meridian
for $K$ commuting with $\lambda$ and we denote by $P$ the
peripheral subgroup of $\pi_1(S^3\setminus K)$ generated by
them.\\[\baselineskip]
\noindent{\bf Proof of Proposition \ref{longnotdead}:}~Suppose
that $K$ is a small hyperbolic knot, $K'$ is a satellite knot, and
that there exists an epimorphism
$$\varphi: \pi_1(S^3\setminus K) \rightarrow  \pi_1(S^3\setminus K').$$

Suppose that $\varphi(\lambda)\neq 1$.
Since a
knot group is torsion-free, $\varphi(P)$ is either
infinite cyclic or isomorphic to ${\bf Z}\oplus{\bf Z}$,
Assume that the former case holds. Then there
is some primitive slope $r=\mu^m\lambda^n$ such that $\varphi(r)=1$,
so $\varphi$ factors through the fundamental group of $K(r)$. This is
impossible, since
by assumption, $r \ne \lambda^{\pm 1}$, so
$\pi_1(K(r))$ has finite abelianization, and thus cannot surject onto $\pi_1(S^3\setminus K')$.\\

Thus we can assume that $\varphi(P)\cong {\bf Z}\oplus{\bf Z}$.
Suppose $f: E(K)\to E(K')$ is a map realizing $\varphi$ and let
$T=\partial E(K)$. Let $T'$ be a JSJ torus of $E(K')$.

By the enclosing property of the JSJ decomposition we may assume
that $f$ has been homotoped so that

(1) $f(T) \subset \Sigma$, where $\Sigma$ is a piece of the JSJ
decomposition.

\noindent Moreover we can assume that

(2) $f^{-1}(T')$ is a 2-sided incompressible surface in $E(K)$; and

(3) $f^{-1}(T')$ has minimum number of components.

Note that $f^{-1}(T')$ can not be empty, otherwise since $T'$ is a separating
torus in $E(K')$, $f(E(K))$ will miss some vertex manifold of
$E(K')$, therefore $f_*=\varphi$ cannot be surjective. No component
$T^*$ of $f^{-1}(T')$ is parallel to $T$, otherwise we can push the
image of the product bounded by $T$ and $T^*$ crossing $T'$ to
reduce the number of components of $f^{-1}(T')$. Therefore
$f^{-1}(T')$ is closed embedded
essential surface in $E(K)$. This is false since $K$ is small.\qed

\subsection{Property L}

We start by proving Proposition \ref{nontrivlong} which shows that
Property L allows control of the image of a longitude under a knot
group epimorphism.
Notation for a longitude and meridian is that of \S 3.3. We remark
that it is here that crucial use is made of \cite{KM}.\\[\baselineskip]
\noindent{\bf Proof of Proposition \ref{nontrivlong}:}~Let $K$ be a
hyperbolic knot and
$\varphi:\pi_1(S^3\setminus K)\rightarrow \pi_1(S^3\setminus K')$ an
epimorphism. If $\varphi(\lambda) = 1$, then the epimorphism $\varphi$
factorizes through an epimorphism $\varphi':\pi_1(K(0))\rightarrow
\pi_1(S^3\setminus K')$. Now Theorem \ref{KM} provides a curve of
characters $C \subset X(K')$ whose generic point is the character of
an irreducible representation. By Lemma \ref{one-to-one}, the curve $D
= {\varphi '}^*(C) \subset X(K(0))$ contradicts the
Property L assumption.\qed\\[\baselineskip]
Together with Proposition \ref{longnotdead}, we obtain:

\begin{corollary}
\label{nosatellite} There cannot be an epimorphsim from the
group of a small knot having Property L onto the
group of a satellite knot.
\end{corollary}

Now we prove Proposition \ref{criterpropertyL} whose content is given by the following Lemmas:

\begin{lemma}
\label{criterparabolic}
Let $K$ be a small hyperbolic knot and $\rho$ any parabolic representation.
Suppose that $\rho(\lambda)\neq 1$, then Property $L$ holds for $K$.
\end{lemma}

\noindent{\bf Proof:}~Suppose that $X(K(0))$ contained a curve $C$ of
characters of irreducible representations. Then the epimorphism
$\psi : \pi_1(S^3\setminus K)\rightarrow \pi_1(K(0))$ induced by
$0$-Dehn surgery together with Lemma \ref{one-to-one}, provides a curve $\psi
^*(C)=D\subset X(K)$. Proposition \ref{prepchars} shows that $D$
contains a p-rep. character $\chi_\theta$, and by assumption
$\theta(\lambda)\neq 1$. On the other hand
$\chi_\theta=\psi^*(\chi_{\theta'})=\chi_{\theta'\psi}$ for
some $\chi_{\theta'}$ in $C$. Hence $\theta=\theta'\psi$ up to
conjugacy. Since $\theta$ factorizes through $\psi:
\pi_1(S^3\setminus K)\rightarrow \pi_1(K(0))$, we must have
$\theta(\lambda)= 1$ and therefore reach a contradiction.\qed\\

The second part of Proposition \ref{criterpropertyL} follows from:

\begin{lemma}
\label{notstrictboundaryslope}
Let $K$ be a small hyperbolic knot.
If the longitude is not a strict boundary slope,   then
Property L  holds for $K$.
\end{lemma}

\noindent{\bf Proof:}~ Let $K$ be a small hyperbolic knot whose
preferred longitude $\lambda$ for $K$ is not a strict boundary
slope. Assume that the character variety $X(K(0))$ contains a curve of
characters $C$ whose generic element is the character of an
irreducible representation. The epimorphism $\varphi:
\pi_1(S^3\setminus K)\rightarrow \pi_1(K(0))$ and Lemma
\ref{one-to-one} provides a curve component $D = \varphi^*(C) \subset
X(K)$.  Since $\varphi(\lambda) = 1$, $f_\lambda: D \to {\bf C}$ is
identically $0$. Thus we
deduce that $\lambda$ is a boundary slope detected by any ideal point
of $D$ (cf. proof of Lemma \ref{move_meridian}).

Fix an irreducible character $\chi_\rho \in D$. By hypothesis, $\lambda$ is not a strict boundary slope, so \cite[Proposition 4.7(2)]{BoZ1} implies that the restriction of $\rho$ to the index $2$ subgroup $\tilde \pi$ of $\pi_1(S^3 \setminus K)$ has Abelian image. The irreducibility of $\rho$ implies that this image is non-central in $\SL(2;{\bf C})$, and as it is normal in the image of $\rho$, the latter is conjugate into the subgroup of $\SL(2;{\bf C})$ of matrices which are either diagonal or have zeroes on the diagonal. Further, the image of $\tilde \pi$ conjugates into the diagonal matrices and that of a meridian of $K$ conjugates to a matrix with zeroes on the diagonal. Any such representation of $\pi_1(S^3 \setminus K)$ has image a finite binary dihedral group. As there are only finitely many such characters of $\pi_1(S^3 \setminus K)$ (\cite[Theorem 10]{Kl}), the generic character in $D$, hence $C$, is reducible, a contradiction.   \qed \\

We can now give the proof of our main technical result Theorem \ref{mainsmall}.\\

\noindent{\bf Proof of Theorem \ref{mainsmall}:}~ We are supposing
that $K$ is a small hyperbolic knot with Property L.
By Corollary \ref{nosatellite}, the targets cannot be fundamental groups of satellite knot complements, hence they must be fundamental groups of hyperbolic or torus knot complements. Then the proof follows from Theorem \ref{smallsimonintro}.
\qed\\[\baselineskip]

\section{Results for 2-bridge knots}

Given the discussion for torus knots in \S 3.1, it suffices to deal
with the case of a hyperbolic 2-bridge knot.

\subsection{Proof of Theorem \ref{main} and Proposition \ref{2bridgehasL}}

As mentioned in the introduction  Theorem \ref{main}
follows from Theorem \ref{mainsmall} and Proposition
\ref{2bridgehasL}. This Proposition is a straightforward consequence of Lemma \ref{criterparabolic} and
of the
following lemma of Riley (see Lemma 1 of \cite{Ri2}).
We have decided to include a proof of this lemma since it is a crucial point.

\begin{lemma}
\label{prepdoesntkilllong}
If  $\theta: \pi_1(S^3\setminus K) \to \PSL(2,{\bf C})$ is a p-rep,
then $\theta(\lambda)\neq 1$.\end{lemma}

\noindent{\bf Proof:}~We begin by recalling some of the basic set up of p-reps. of 2-bridge knot
groups (see \cite{Ri1}).  Let $K$ be 2-bridge of normal form $(p,q)$, so $p$ and $q$ are odd integers such that $0 < q < p$.  The case of $q = 1$ is that of 2-bridge torus knots. The group $\pi_1(S^3\setminus K)$ has a presentation

$$<x_1,x_2~|~wx_1w^{-1} = x_2^{-1}>,~~\hbox{where}~~x_1, x_2~~\hbox{are meridians and}$$

\noindent $w=w(x_1,x_2)$ is given by $x_1^{\epsilon_1}x_2^{\epsilon_2}\ldots x_1^{\epsilon_{p-2}}x_2^{\epsilon_{p-1}}$.
Furthermore, each exponent $\epsilon_j = (-1)^{[jq/p]}$ where $[x]$ denotes the integer part of $x$, and
$\epsilon_j=\epsilon_{p-j}$. Hence $\sigma=\Sigma \epsilon_j$ is even.

The standard form for a p-rep sends the meridians $x_1$ and $x_2$ to parabolic elements:

$$ \pmatrix{1 &  1\cr 0 & 1}~~\hbox{and}~~ \pmatrix{1 &  0\cr -y & 1},$$
for some non-zero algebraic integer $y$ (indeed $y$ is a unit). The relation in the presentation provides
a p-rep polynomial $\Lambda(y)$, and all p-reps  determine and are
determined by solutions to $\Lambda(y)=0$.  The image of $w$ under p-rep has the form
$$W = \pmatrix{0 &  w_{12}\cr w_{21} & w_{22}},$$
with the entries being functions of the variable $y$. In addition, as is shown in \cite{Ri1}, the image
of a longitude that commutes with $x_1$ has the form
$$\pmatrix{1 &  -2g\cr 0 & 1}~~\hbox{for some algebraic integer}~~g=g(y).$$
Indeed, as shown in \cite{Ri1}, $g=w_{12}w_{22}+\sigma$.  Thus, to prove the lemma we need to show that
$g=g(y)\neq 0$.

This is done as follows. First, observe that (mod $2$), the matrix $W$
for the 2-bridge knot of normal form $(p,q)$ is the same as the matrix
$W'$ one obtains from the 2-bridge torus knot with normal form
$(p,1)$.  Furthermore, the word $w$ in the case of $(p,1)$ is given as
$(x_1x_2)^n$ with $n=(p-1)/2$ the degree of $\Lambda(y)$. Using this
allows for an easy recursive definition of the matrix $W'$ in this
case (see \S 5 of \cite{Ri1}); namely define two sequences of
polynomials $f_j=f_j(y)$ and $g_j=g_j(y)$, with $f_0(y)=g_0(y)=1$ and:

$$f_{j+1}(y) = f_j(y) + yg_j(y)~~\hbox{and}~~g_{j+1} = f_{j+1} + g_j(y)$$
Then the matrix $W'$ is given by:

$$W' = \pmatrix{f_n &  g_{n-1}\cr yg_{n-1} & f_{n-1}}.$$

\noindent In particular, the p-rep condition implies $f_n(y)=0$. Using the recursive formula, we have
$f_n(y) = f_{n-1}+yg_{n-1}(y)$, and
the p-rep condition (i.e $f_n(y)=0$) means that the matrix $W'$ is given by

$$W' = \pmatrix{0 &  g_{n-1}\cr yg_{n-1} & -yg_{n-1}},$$

\noindent which (mod $2$) is

$$W' = \pmatrix{0 &  g_{n-1}\cr yg_{n-1} & yg_{n-1}}.$$

\noindent We deduce from these comments that $w_{12}w_{22} = -w_{12}w_{21}$ (mod $2$).  The latter is 1
since it is the determinant of $W$. As noted above,  $\sigma$ is even, hence, it follows that $g=w_{12}w_{22}+\sigma$ is congruent to 1 (mod $2$), and so in particular is not zero as required.\qed

\subsection{Proof of Corollary \ref{nonzerodegree}}

\noindent{\bf Proof of Corollary \ref{nonzerodegree}:}~As before we
let $\mu$ and $\lambda$ denote a meridian and a longitude of $K$.
Firstly, we note that if $K$ is a hyperbolic 2-bridge knot and
$\varphi:\pi_1(S^3\setminus K)\rightarrow \pi_1(S^3\setminus K')$ is
an epimorphism, then Proposition \ref{longnotdead} and Lemma
\ref{prepdoesntkilllong} combine to show that $K'$ is either a
hyperbolic or torus knot.

In the case of $K'$ a hyperbolic knot, since $\varphi(\lambda)\neq 1$,
the epimorphism is non-degenerate in the sense of \cite{BRW}, and in
particular $\varphi(\mu)$ is a peripheral element of $\pi_1(S^3\setminus K')$. Hence,  \cite[Theorem 3.15]{BB}
applies to show that $K'$ is also a 2-bridge knot.  Furthermore, as noted in the proof of Corollary 6.5 of
\cite{BRW}, the homomoprhism $\varphi$ is induced by a map of non-zero degree.

In the case when $K'$ is a torus knot, $\varphi(\lambda)$ (and
therefore also $\varphi(\mu)$), need not be a peripheral element.
Suppose that $K'$ is an $(r,s)$-torus knot and fix a meridian $\mu'$
of $K'$.  There is a homomorphism $\psi:\pi_1(S^3\setminus
K')\rightarrow C_{r,s}$ (where as in \S 3.1, $C_{r,s}$ denotes the
free product of two cyclic groups of orders $r$ and $s$) and
generators $a$ of ${\bf Z}/r$ and $b$ of ${\bf Z}/s$ for which
$\psi(\mu') = ab$. Theorem 2.1 of \cite{GAR1} and the remark following
it, shows that one of $r$ or $s$ equals $2$, say $r=2$.  In particular
$K'$ is a 2-bridge torus knot.  We finish off this case as we did the
previous one using \cite{BRW} once we show that
$\varphi(<\mu,\lambda>)$ is a subgroup of finite index in the
peripheral subgroup $P'$ of $\pi_1(S^3\setminus K')$.  To do this, it
suffices to show $\varphi(\mu)$ is a meridian of $K'$ since the
centraliser of $\mu'$ in $\pi_1(S^3 \setminus K')$ is $P'$.

To that end, Theorem 1.2 of \cite{GAR1} shows that there is an
isomorphism $\theta: C_{2,s} \to C_{2,s}$ such that $\theta \psi
\varphi(\mu) = a b^m$ for some integer $m$. Up to inner isomorphism,
we can suppose that $\theta(a) = a$ and $\theta(b) = b^k$ for some
$k$ coprime with $s$ (see for example \cite[Theorem 13(1), Corollary
14]{GAR2}). Thus we can assume that $\psi \varphi(\mu) = ab^m$. Now
$\varphi(\mu)$ equals $(\mu')^{\pm 1}$ up to multiplication by a
commutator, so abelianizing in $C_{2,s}$ shows that $m \equiv \pm 1$
(mod $s$). Hence $\psi \varphi(\mu) = ab^{\pm 1}$, so up to
conjugation in $C_{2,s}$, $\psi \varphi(\mu) = (ab)^{\pm 1} =
\psi(\mu')^{\pm 1}$. It follows that $\varphi(\mu) = h^k (\mu')^{\pm
1}$ where $h \in \pi_1(S^3 \setminus K')$ is the fibre class. Since
$K'$ is non-trivial, $|s| \geq 3$, and so as $h$ represents $2s$ in
$H_1(S^3 \setminus K') \cong {\bf Z}$, it must be that $k = 0$. Thus
$\varphi(\mu) = (\mu')^{\pm 1}$, which completes the proof.\qed\\[\baselineskip]
We conclude this section with some remarks on the proof of a stronger
version of Corollary \ref{nonzerodegree}.  Before stating this
result, we recall that if $G$ and $H$ are groups and
$\varphi :G\rightarrow H$ is a homomorphism, then $\varphi$ is called
a {\em virtual epimorphism} if $\varphi(G)$ has finite index in $H$.

\begin{theorem}
\label{nonzerodegree2} Let $K$ be a 2-bridge hyperbolic knot, $K'$ be
a non-trivial knot. If there is a virtual epimorphism
$\varphi:\pi_1(S^3\setminus K)\rightarrow \pi_1(S^3\setminus K')$,
then $\varphi$ is induced by a map  $f: S^3\setminus K \rightarrow
S^3\setminus K'$ of non-zero degree. Furthermore, $K'$ is
necessarily a 2-bridge knot, and $\varphi$ is surjective if $K'$ is
hyperbolic.
\end{theorem}

\noindent {\bf Sketch of the Proof:}~Since a subgroup of finite index
in a satellite knot group continues to contain an essential
${\bf Z}\oplus {\bf Z}$ the proof of Proposition \ref{longnotdead}
can be applied to rule out the case of satellite knot groups as
targets.

In the case where the targets are hyperbolic, we can deduce
that this virtual epimorphism is an epimorphism and
we argue as before; briefly, since the
peripheral subgroup is mapped to a ${\bf Z}\oplus {\bf Z}$ in the image,
and since $K$ is 2-bridge, it follows from \cite[Corollary 5]{BW} that
these image groups are 2-bridge knot groups (being generated
by two conjugate peripheral elements). However, it is well-known
that a 2-bridge hyperbolic knot complement has no free symmetries,
and so cannot properly cover any other
hyperbolic 3-manifold (see \cite{Sa} for example).

When the targets are torus knot groups, standard considerations show
that the image
of $\varphi$ is the fundamental group of a Seifert Fiber Space
with base orbifold a disc with cone points. Moreover, this 2-orbifold
group is generated by the images of the two conjugate meridians of $K$.
It is easily seen that this forces the base orbifold to be a disc with
two cone points. It now follows
from \cite[Proposition 17]{GAR2} that the
base orbifold group is $C_{2,s}$ where
$s$ is odd, and the proof is completed as before.\qed\\[\baselineskip]
\noindent{\bf Remark:}~Note that the paper
\cite{ORS} gives a systematic construction of epimorphisms between
2-bridge knot groups. In particular the epimorphisms constructed
by the methods of \cite{ORS} are induced by maps of non-zero degree.
Corollary \ref{nonzerodegree} and Theorem \ref{nonzerodegree2}
show that in fact {\em any}
(virtual) epimorphism from a 2-bridge knot group to any knot group
is induced by a map of non-zero degree.

\subsection{Minimal manifolds and Simon's Conjecture}~The methods of
this paper also prove the following strong form of Conjecture \ref{simons} in
certain cases.

\begin{theorem}
\label{onecpt}
Suppose $K\subset S^3$ is a hyperbolic knot for which the canonical component of $X(K)$ is the
only component that contains the character of an irreducible representation.  Then
Conjecture \ref{simons} holds for $K$.\end{theorem}

This obviously follows from the stronger theorem stated below.

\begin{theorem}
\label{onecptII}
Suppose $K\subset S^3$ is as in Theorem \ref{onecpt}.  Then
$\pi_1(S^3\setminus K)$ does not surject onto the fundamental group of any other non-trivial knot complement.
\end{theorem}

\noindent{\bf Proof:}~Assume to the contrary that
$\phi:\pi_1(S^3\setminus K)\rightarrow \pi_1(S^3\setminus K')$ is a
surjection. It will be convenient to make use of the $\PSL(2,{\bf C})$
character variety. Let $Y_0(K)$ denote the canonical component of
$Y(K)$.

$K'$ cannot be a torus knot since $Y_0(K)$ contains the character of a
faithful representation of $\pi_1(S^3\setminus K)$ and $Y(C_{p,q})$
clearly contains no such character.  That is to say $Y_0(K)\neq
\phi^*(Y(C_{p,q}))$.

Theorem \ref{numberofcptsbdd} handles the case when $K'$ is
hyperbolic. More precisely, taking $G=\pi_1(S^3\setminus K)$, in the
notation of Theorem \ref{numberofcptsbdd}, $k\leq 1$. Since, $G$
surjects onto itself, we deduce that there can be no other knot group
quotient.

Now assume that $K'$ is a satellite knot.  In this case, we use
Theorem \ref{KM} to deduce that $\phi^*(Y(K'))$ coincides with $Y(K)$.
However, if $\chi_\rho$ denotes the character of the faithful discrete
representation on the canonical component $Y_0(K)$,
then there is a character $\chi_\nu \in Y(K')$ with $\rho = \nu\phi$. But this is clearly impossible.\qed\\[\baselineskip]
By \cite{Mat}, when $n\geq 7$ is not divisble by $3$, the
$(-2,3,n)$-pretzel knot satisfies the hypothesis of Theorem
\ref{onecpt}. Hence we get.

\begin{corollary}
\label{pretzelsimon}
Suppose that $n\geq 7$ is not divisble by $3$, then Conjecture \ref{simons} holds for the $(-2,3,n)$-pretzel knot.
\end{corollary}

\section{Possible extension of Simon's Conjecture}

We first state a possible extension of Simon's Conjecture for links.
To that end, recall that a
{\em boundary link} is a link whose components bound
disjoint Seifert surfaces.  Such a link (say with n components) has
a fundamental group that surjects onto a non-abelian free group of
rank $n$.  A {\em homology boundary link} of $n$ components is a
link of $n$ components
whose fundamental group surjects onto a non-abelian free group of rank $n$.

\begin{conjecture}
  \label{nSimon} Let $L\subset S^3$ be a non-trivial link of $n \geq
  2$ components. If $\pi_1(S^3\setminus L)$ surjects onto infinitely
  many distinct link groups of $n$ components, then $L$ is a homology
  boundary link.\end{conjecture}

This conjecture is motivated by Simon's conjecture for knots and the
following observations: \\[\baselineskip]
\noindent If $n\geq 2$, then the trivial link of
$n$-components has a fundamental group which is free of rank $n\geq
2$. Hence, it surjects onto all link groups which are generated by $n$
elements. This argument can now be made by replacing the trivial link
by a homology boundary link. In particular, since there are
non-trivial boundary links of $2$ components, the fundamental groups
of such link complements will surject onto all two components 2-bridge
link groups.

Hyperbolic examples are easily constructed from this
using \cite{Ka} for example. Hence the group of any link with $n$ components
is the homomorphic image of the fundamental group of a hyperbolic link
with $n$ components.\\[\baselineskip]
As in Corollary \ref{highdimcpt}, Theorem \ref{numberofcptsbdd}
provides information about the dimension of the character variety of a
homology boundary link with $n \geq 2$ components (see also \cite{CL}):

\begin{corollary}
$\dim (X(L)) > n$  for each homology boundary link  $L$ of $n\geq 2$ components.
\end{corollary}

\noindent{\bf Proof:}~Let $L$ be a homology boundary link of $n$
components. The group of $L$ surjects onto all $n$ component
link groups  which are generated by $n$ elements. As we note in
the Remark following this proof, infinitely
many of these correspond to distinct hyperbolic link complements, and so
$\dim (X(L)) \ge n$ by Lemma \ref{one-to-one}. Hence $\dim (X(L)) > n$ by
Theorem \ref{numberofcptsbdd}.\qed\\[\baselineskip]
\noindent{\bf Remark:}~It is easy to see that there are infinitely
many $n$-component hyperbolic links whose groups are generated by
$n$-elements. Briefly, by Thurston's hyperbolisation theorem for surface bundles
(\cite{Ot}, \cite{Th}) a pseudo-Anosov pure braid with $n-1$ strings
together with its axis forms a hyperbolic $n$-bridge link with $n$
components. Moreover the group of such a link is generated
by $n$ elements. Since there are infinitely many conjugacy classes of
pseudo-Anosov pure braids with $n-1$ strings, infinitely many distinct
hyperbolic link complements can be obtained in this way.\qed \\[\baselineskip]
Another natural extension of Simon's Conjecture is.

\begin{conjecture}
\label{simon_general} Let $X$ be a knot exterior in a closed
orientable 3-manifold for which $H_1(X:{\bf Q})\cong {\bf Q}$.  Then
$\pi_1(X)$ surjects onto only finitely many groups $\pi_1(X_i)$
where $X_i$ is a knot exterior with $H_1(X_i:{\bf Q})\cong {\bf
Q}$.\end{conjecture}

The condition on the rational homology is clearly a necessary
condition (otherwise one can use surjections that factor through a
non-abelian free group once again). Even here little seems known.
Indeed, even for small manifolds as in Conjecture
\ref{simon_general} we cannot make as much progress as in the case
of $S^3$, since Theorem \ref{KM} of Kronheimer and Mrowka is not
known to hold in this generality.

\bigskip
 \noindent Laboratoire \'Emile Picard, CNRS UMR 5580,\\
 Universit\'e Paul Sabatier,\\
 F-31062 Toulouse Cedex 4, France.

\noindent Email:~boileau@picard.ups-tlse.fr\\[\baselineskip]
D\'epartment de Mathematiques,\\
U.Q.A.M.\\
P. O. Box 8888, Centre-ville,\\
Montr\'eal, Qc, H3C 3P8, Canada.

\noindent Email:~boyer@math.uqam.ca\\[\baselineskip]
Department of Mathematics,\\
University of Texas\\
Austin, TX 78712, USA.

\noindent Email:~areid@math.utexas.edu\\[\baselineskip]
LAMA Department of Mathematics,\\
Peking University,\\
Beijing 100871, China.

\noindent Email:~wangsc@math.pku.edu.cn\\

\end{document}